\documentclass[a4paper,12pt]{article}
\usepackage{latexsym}
\usepackage{amsfonts}
\usepackage{xypic}
\newtheorem{theorem}{Theorem}[section]
\newtheorem{corollary}[theorem]{Corollary}
\newtheorem{definition}{Definition}
\newtheorem{example}[theorem]{Example}
\newtheorem{remark}[theorem]{Remark}
\newtheorem{lemma}[theorem]{Lemma}
\newtheorem{proposition}[theorem]{Proposition}
\title{On p-closed spaces\thanks{1991 Math.\ Subject
Classification --- Primary: 54D20, 54D25, Secondary: 54A05,
54D30, 54H05. \protect\newline Key words and phrases --- p-closed,
QHC, strongly compact, nearly compact, preopen, consolidation.
\protect\newline Research supported partially by the Ella and Georg
Ehrnrooth Foundation at Merita Bank, Finland and by the
Japan - Scandinavia Sasakawa Foundation.}}
\author{Julian {\sc Dontchev}, Maximilian {\sc Ganster} and
Takashi {\sc Noiri}}
\date{}
\begin{document}
\baselineskip=20pt plus 1pt minus 1pt
\newcommand{\fxy}{$f \colon (X,\tau) \rightarrow (Y,\sigma)$}
\newcommand{\dpo}{$\delta$-preopen}
\newcommand{\dpc}{$\delta$p-closed}
\newcommand{\dpco}{$\delta$p-compact}
\newcommand{\ptc}{pre-$\theta$-converge}
\newcommand{\pta}{pre-$\theta$-accumulate}
\maketitle
\begin{abstract}
In this paper we will continue the study of p-closed spaces. This
class of spaces is strictly placed between the class of strongly
compact spaces and the class of quasi-H-closed spaces. We will
provide new characterizations of p-closed spaces and investigate
their relationships with some other classes of
topological spaces.
\end{abstract}

\section{Introduction and Preliminaries}\label{s1}

The aim of this paper is to continue the study of p-closed
spaces, which were introduced in 1989 by Abo-Khadra \cite{A1}.
A topological space $(X,\tau)$ is called {\em p-closed} if every
preopen cover of $X$ has a finite subfamily whose pre-closures
cover $X$.

Let $A$ be a subset of a topological space $(X,\tau)$. Following
Kronheimer \cite{K1}, we call the interior of the closure of $A$,
denoted by $A^+$, the {\em consolidation} of $A$. Sets included
in their consolidation play a significant role in e.g.\ questions
concerning covering properties, decompositions of continuity,
etc. Such sets are called {\em preopen} \cite{MAH1} or {\em
locally dense} \cite{CM1}. A subset $A$ of a space $(X,\tau)$ is
called {\em preclosed} if its complement is preopen, i.e. if
${\rm cl}({\rm int} A) \subseteq A$. The preclosure of $A
\subseteq X$, denoted by ${\rm pcl} (A)$, is the intersection of
all preclosed supersets of $A$. Since any union of preopen sets
is also preopen, the preclosure of every set is preclosed. It is
well known that ${\rm pcl} A = A \cup {\rm cl}({\rm int} A)$ for
any $A \subseteq X$.

Another interesting property of preopen sets is the following:
When a certain topological property is inherited by both open and
dense sets, it is often then inherited by preopen sets.

Several important concepts in Topology are and can be defined in
terms of preopen sets. Among the most well-known are Bourbaki's
submaximal spaces (see \cite{AC1}). A topological space is called
{\em submaximal} if every (locally) dense subset is open or,
equivalently, if every subset is locally closed,i.e. the
intersection of an open set and a closed set. Another class of
spaces commonly characterized in terms of preopen sets is the
class of strongly irresolvable spaces introduced by Foran and
Liebnitz in \cite{FL1}. A topological space $(X,\tau)$ is called
{\em strongly irresolvable} \cite{FL1} if every open subspace of
$X$ is irresolvable, i.e.\ it cannot be represented as the
disjoint union of two dense subsets. Subspaces that contain
two disjoint dense subsets are called {\em resolvable}. Ganster
\cite{GA1} has pointed out that a space is strongly irresolvable
if and only if every preopen set is semi-open, where a subset $S$
of a space $(X,\tau)$ is called {\em semi-open} if \ $S \subseteq 
{\rm cl}({\rm int} S)$. We will denote the families preopen
(resp.\ semi-open) sets of a space $(X,\tau)$ by $PO(X)$ (resp.\
$SO(X)$).

Many classical topological notions such as compactness and
connectedness have been extended by using preopen sets instead
of open sets. Among them are the class of {\em strongly compact
spaces} \cite{MM1} (= every preopen cover has a finite subcover)
studied by Jankovi\'{c}, Reilly and Vamanamurthy \cite{JRV1} and
by Ganster \cite{GA2}, and the class of {\em preconnected spaces}
(= spaces that cannot be represented as the disjoint union of two
preopen subsets) introduced by Popa \cite{P1}. The study of
topological properties via preopenness has gained significant
importance in General Topology and one example for that is the
fact four (out of the ten) articles in the 1998 Volume of
"Memoirs of the Faculty of Science Kochi University Series A
Mathematics" were more or less devoted to preopen sets.

A point $x \in X$ is called a {\em $\delta$-cluster point} of a
set $A$ \cite{V1} if $A \cap U \not= \emptyset$ for every regular
open set $U$ containing $x$. The set of all $\delta$-cluster
points of $A$ forms the {\em $\delta$-closure} of $A$ denoted by
${\rm cl}_{\delta} (A)$, and $A$ is called {\em $\delta$-closed}
\cite{V1} if $A = {\rm cl}_{\delta} (A)$. If $A \subseteq {\rm
int} ({\rm cl}_{\delta} (A))$, then $A$ is said to be {\em
$\delta$-preopen} \cite{RM2}. Complements of $\delta$-preopen
sets are called {\em $\delta$-preclosed} and the
{\em $\delta$-preclosure} of a set $A$, denoted by $\delta$-${\rm
pcl} (A)$, is the intersection of all $\delta$-preclosed
supersets of $A$.

Following \cite{RM1}, we will call a topological space $(X,\tau)$
{\em $\delta$p-closed} if for every $\delta$-preopen cover
$\{V_{\alpha} \colon \alpha \in A \}$ of $X$, there exists a
finite subset $A_0$ of $A$ such that $X = \cup \{ {\delta}$-${\rm
pcl}(V_{\alpha}) \colon \alpha \in A_0 \}$.

\section{p-closed spaces}\label{s2}

\begin{definition}\label{d1}
{\em A topological space $(X,\tau)$ is said to be {\em p-closed}
\cite{A1} (resp.\ {\em quasi-H-closed} = QHC) if for every
preopen (resp.\ open) cover $\{ V_{\alpha} \colon \alpha \in A
\}$ of $X$, there exists a finite subset $A_0$ of $A$ such that
$X = \cup \{ {\rm pcl}(V_{\alpha}) \colon \alpha \in A_0 \}$
(resp. $X = \cup \{ {\rm cl}(V_{\alpha}) \colon \alpha \in A_0
\}$}.
\end{definition}

It is clear that every strongly compact space is p-closed, and
that every p-closed space is QHC. We also observe that a space
$(X,\tau)$ is QHC if and only if every preopen cover has a finite
dense subsystem (= finite subfamily whose union is a dense
subset). Since every preopen set is $\delta$-preopen, we have
${\delta}$-${\rm pcl} S \subseteq {\rm pcl} S$ for every $S
\subseteq X$ . This implies that every  $\delta$p-closed space
is p-closed.

\begin{theorem}\label{t1}
Let $(X,\tau)$ be QHC and strongly irresolvable. Then $(X,\tau)$
is p-closed.
\end{theorem}

{\em Proof.}  Let $\{ V_i  \colon i \in I \}$ be any preopen
cover of $X$. Since $X$ is QHC, there exists a finite subset $A$
of $I$ such that $X = \cup \{ {\rm cl} (S_i) : i \in A \}$. Since
$X$ is strongly irresolvable, $S_i \in SO(X)$ and therefore ${\rm
cl} (S_i) = {\rm cl}({\rm int} (S_i)) = {\rm pcl} (S_i)$ for each
$i \in I$. Hence $X$ is p-closed. $\Box$

\begin{corollary}\label{c1}
Let  $(X,\tau)$ be strongly irresolvable (or submaximal). Then 
$(X,\tau)$ is p-closed if and only if it is QHC .
\end{corollary}

Observe that a p-closed space need not be strongly irresolvable
as any finite indiscrete space shows. However, we do have the
following result.

\begin{theorem}\label{t2}
Let  $(X,\tau)$ be a p-closed $T_0$ space. Then $(X,\tau)$ is
strongly irresolvable.
\end{theorem}

{\em Proof.} Suppose that $W$ is a nonempty, open and resolvable
subspace of $X$. Then $W$ is dense-in-itself and also infinite,
since $(X,\tau)$ is $T_0$. Let $W = E_1 \cup E_2$, where $E_1$
and $E_2$ are disjoint dense subsets of $W$, and wlog.\ we may
assume that $E_1$ is infinite. Moreover, let $A = \{ x \in E_1
\colon \{ x \} \in PO(X) \} $. Observe that for each $y \in E_1
\setminus A$ , $\{ y \}$ is nowhere dense. Now pick $y \in E_1
\setminus A$. If $S_y = (X \setminus W) \cup E_2 \cup \{ y \}$
then $S_y$ is dense and therefore preopen. If $G$ is a nonempty
open set contained in $S_y$, then $G \cap E_1 \subseteq \{ y \}$
and so $G \cap W \subseteq {\rm cl} (E_1) \subseteq {\rm cl} \{
y \}$. Since $\{ y \}$ is nowhere dense, $G \cap W$ is empty and
so ${\rm cl}({\rm int} (S_y)) \subseteq X \setminus W$, thus $pcl
S_y = S_y$. Now, observe that $\{ \{ x \} \colon x \in A \} \cup
\{ S_y \colon y \in E_1 \setminus A \}$ is a preopen cover of
$X$. Hence there exists a finite subset $A_1$ of $A$ and a finite
subset $A_2$ of $E_1 \setminus A$ such that $X = \{ \{ x \}
\colon x \in A_1 \} \cup \{ S_y \colon y \in A_2 \}$. Then, $E_1
\subseteq A_1 \cup A_2$ which is a contradiction. Thus $X$ is
strongly irresolvable. $\Box$

\medskip
By combining the previous two results we immediately have:

\begin{theorem}\label{t3}
Let $(X,\tau)$ be a $T_0$ space. Then $(X,\tau)$ is p-closed if
and only if $(X,\tau)$ is QHC and strongly irresolvable.
\end{theorem}

The following diagram exhibits the relationships between the
class of p-closed spaces and some related classes of topological
spaces. Note that none of the implications is reversible.

$$
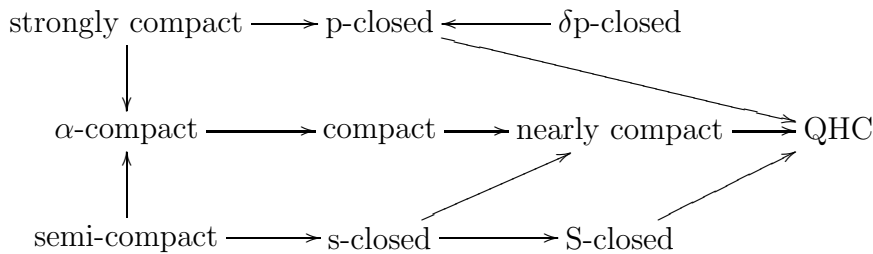
\diagram
\text{strongly compact} \rto \dto & \text{p-closed} \drrto &
\text{\dpc} \lto \\ \text{$\alpha$-compact} \rto &
\text{compact} \rto & \text{nearly compact} \rto & \text{QHC}
\\ \text{semi-compact} \uto \rto & \text{s-closed} \rto \urto
& \text{S-closed} \urto
\enddiagram
$$

\begin{example}\label{e1}
{\em (i) Recall that a space $(X,\tau)$ is called {\em
$\alpha$-scattered} \cite{DGR1} if it has a dense set of isolated
points. Clearly every $\alpha$-scattered space is strongly
irresolvable and so, by Theorem~\ref{t1}, every
$\alpha$-scattered QHC space is p-closed. In particular, the
Katetov extension $\kappa {\mathbb N}$ of the set of natural
numbers ${\mathbb N}$ (see e.g.\ \cite{PW}) is p-closed and not
compact, hence not strongly compact.

(ii) The unit interval $[0,1]$ with the usual topology is
compact, hence QHC, but not p-closed since it is resolvable.

(iii) Let $X = {\mathbb R}$, $\tau = \{ \emptyset, \{ 0 \}, X
\}$. Then, $X$ is p-closed and s-closed but not $\alpha$-compact
and hence not strongly compact (A space is {\em $\alpha$-compact}
if every cover by $\alpha$-open sets has a finite subcover, where
a set is {\em $\alpha$-open} if it is the difference of an open
and a nowhere dense set; clearly every $\alpha$-open set is
preopen but not vice versa). Additionally, this space is not
$\delta$p-closed since every subset is \dpo.}
\end{example}

We next discuss the relationship between p-closedness and
compactness. Recall that a space $(X,\tau)$ is called {\em nearly
compact} \cite{SM1} if every cover of $X$ by regular open sets
has a finite subcover, i.e.\ the semiregularization $(X,\tau_s)$
of $(X,\tau)$ is compact. Example 4.8 (d) in \cite{PW} shows that
there exists a Hausdorff, non-compact, semi-regular and QHC space
with a dense set of isolated points. Such a space is p-closed but
not nearly compact. Example 2.10 in \cite{RM1} provides another
such example.

For any infinite cardinal $\kappa$, a topological space
$(X,\tau)$ is called {\em $\kappa$-extremally disconnected} (=
$\kappa$-e.d.) \cite{DG1} if the boundary of every regular open
set has cardinality (strictly) less than $\kappa$. Several
topological spaces share this property for $\kappa = \aleph_{0}$.
Since there are finite spaces which fail to be extremally
disconnected, clearly $\aleph_0$-extremal disconnectedness
is a strictly weaker property than extremal disconnectedness.

\begin{theorem}\label{t4}
If a topological space $(X,\tau)$ is p-closed and
$\aleph_0$-extremally disconnected (resp.\ extremally
disconnected), then $(X,\tau)$ is nearly compact (resp.\
$s$-closed).
\end{theorem}

{\em Proof.} We first prove the case when the space is
$\aleph_0$-extremally disconnected. Let $\{ V_i \colon i \in I
\}$ be any regular open cover of $X$. For each $i \in I$, we have
${\rm pcl} (V_{i}) = V_{i} \cup V_i \cup {\rm cl} ({\rm
int}(V_{i})) = {\rm cl} (V_{i})$. Since $X$ is p-closed, then
there exists a finite $F \subseteq I$ such that $X = \cup_{i \in
F} {\rm cl} (A_{i})$. Note that for each $A_i$, we have ${\rm cl}
(A_i) = B_i \cup C_i$, where $B_i = {\rm int} ({\rm cl}(A_{i}))$
and $C_i = {\rm cl} (A_{i}) \setminus {\rm int} ({\rm
cl}(A_{i}))$. Since $X$ is $\aleph_0$-extremally disconnected,
then $C_i$ is finite for each $i \in F$. Since $B_i = A_i$, for
each $i \in F$, then $\cup_{i \in F} A_i$ covers $X$ but a finite
amount. Hence, $X$ is nearly compact. The proof of the second
part of the theorem is similar to the first one and hence
omitted. $\Box$

On the other hand (see e.g.\ \cite{PW}, page 450) there exist
dense-in-itself, compact and extremally disconnected Hausdorff
spaces. Such spaces are resolvable and hence cannot be p-closed.

A filter base $\cal F$ on a topological space $(X,\tau)$ is said
to {\em \ptc\ } to a point $x \in X$ if for each $V \in PO(X,x)$,
there exists $F \in {\cal F}$ such that $F \subseteq {\rm pcl}
(V)$. A filter base $\cal F$ is said to {\em \pta\ } at $x \in
X$ if ${\rm pcl} (V) \cap F \not= \emptyset$ for every $V \in
PO(X,x)$ and every $F \in {\cal F}$. The {\em preinterior} of a
set $A$, denoted by ${\rm pint} (A)$, is the union of all preopen
subsets of $A$.

\begin{theorem}\label{t41}
For a topological space $(X,\tau)$ the following conditions are
equivalent:

{\rm (a)} $(X,\tau)$ is p-closed,

{\rm (b)} every maximal filter base \ptc s to some point of $X$,

{\rm (c)} every filter base \pta s at some point of $X$,

{\rm (d)} for every family $\{ V_{\alpha} \colon \alpha \in A \}$
of preclosed subsets such that $\cap \{ V_{\alpha} \colon \alpha
\in A \} = \emptyset$, there exists a finite subset $A_0$ of $A$
such that $\cap \{ {\rm pint} (V_{\alpha}) \colon \alpha \in A_0
\} = \emptyset$.
\end{theorem}

{\em Proof.} (a) $\Rightarrow$ (b): Let $\cal F$ be a maximal
filer base on $X$. Suppose that $\cal F$ does not \ptc\ to any
point of $X$. Since $\cal F$ is maximal, $\cal F$ does not \pta\
at any point of $X$. For each $x \in X$, there exists $F_x \in
{\cal F}$ and $V_x \in PO(X,x)$ such that ${\rm pcl} (V_{x}) \cap
F_x = \emptyset$. The family $\{ V_x \colon x \in X \}$ is a
cover of $X$ by preopen sets of $X$. By (a), there exists a
finite number of points $x_1, x_2, \ldots, x_n$ of $X$ such that
$X = \cup \{ {\rm pcl} (V_{x_{i}}) \colon i = 1,2, \ldots, n \}$.
Since $\cal F$ is a filter base on $X$, there exists $F_0 \in
{\cal F}$ such that $F_0 \subseteq \cap \{ F_{x_{i}} \colon i =
1,2, \ldots, n \}$. Therefore, we obtain $F_0 = \emptyset$. This
is a contradiction.

(b) $\Rightarrow$ (c): Let $\cal F$ be any filter base on $X$.
Then, there exists a maximal filter base ${\cal F}_0$ such that
${\cal F} \subseteq {\cal F}_0$. By (b), ${\cal F}_0$ \ptc s to
some point $x \in X$. For every $F \in {\cal F}$ and every $V \in
PO(X,x)$, there exists $F_0 \in {\cal F}_0$ such that $F_0
\subseteq {\rm pcl} (V)$; hence $\emptyset \not= F_0 \cap F
\subseteq {\rm pcl} (V) \cap F$. This shows that $\cal F$ \pta
s at $x$.

(c) $\Rightarrow$ (d): Let $\{ V_{\alpha} \colon \alpha \in A \}$
be any family of preclosed subsets of $X$ such that $\cap \{
V_{\alpha} \colon \alpha \in A \} = \emptyset$. Let $\Gamma(A)$
denote the ideal of all finite subsets of $A$. Assume that $\cap
\{ {\rm pint} (V_{\alpha}) \colon \alpha \in \gamma \} \not=
\emptyset$ for every $\gamma \in \Gamma(A)$. Then, the family
${\cal F} = \{ \cap_{\alpha \in \gamma} {\rm pint} (V_{\alpha})
\colon \gamma \in \Gamma(A) \}$ is a filter base on $X$. By (c),
$\cal F$ \pta s at some point $x \in X$. Since $X \setminus
V_{\alpha} \colon \alpha \in A \}$ is a cover of $X$, $x \in X
\setminus V_{\alpha_{0}}$ for some $\alpha_0 \in A$. Therefore,
we obtain $X \setminus V_{\alpha_{0}} \in PO(x,x)$, ${\rm pint}
(V_{\alpha_{0}} \in {\cal F}$ and ${\rm pcl} (X \setminus
V_{\alpha_{0}}) \cap {\rm pint} (V_{\alpha_{0}}) = \emptyset$.
This is a contradiction.

(d) $\Rightarrow$ (a): Let $\{ V_{\alpha} \colon \alpha \in A \}$
be a cover of $X$ by preopen sets of $X$. Then $\{ X \setminus
V_{\alpha} \colon \alpha \in A \}$ is a family of preclosed
subsets of $X$ such that $\cap \{ X \setminus V_{\alpha} \colon
\alpha \in A \} = \emptyset$. By (d), there exists a finite
subset $A_0$ of $A$ such that $\cap \{ {\rm pint} (X \setminus
V_{\alpha}) \colon \alpha \in A_{0} \} = \emptyset$; hence $X =
\cup \{ {\rm pcl} (V_{\alpha}) \colon \alpha \in A_{0} \}$. This
shows that $X$ is p-closed. $\Box$

\begin{definition}
{\em A topological space $(X,\tau)$ is said to be {\em strongly
p-regular} (resp.\ {\em p-regular} \cite{ElD1}, {\em almost
p-regular} \cite{MN1}) if for each point $x \in X$ and each
preclosed set (resp.\ closed set, regular closed set) $F$ such
that $x \not\in F$, there exist disjoint preopen sets $U$ and $V$
such that $x \in U$ and $F \subseteq V$.}
\end{definition}

\begin{theorem}\label{t42}
If a topological space $X$ is p-closed and strongly p-regular
(resp.\ p-regular, almost p-regular), then $X$ is strongly
compact (resp.\ compact, nearly compact).
\end{theorem}

{\em Proof.} We prove only the case of p-regular spaces. Let $X$
be a p-closed and p-regular space. Let $\{ V_{\alpha} \colon
\alpha \in A \}$ be any open cover of $X$. For each $x \in X$,
there exists an $\alpha(x) \in A$ such that $x \in
V_{\alpha(x)}$. Since $X$ is p-regular, there exists $U(x) \in
PO(X)$ such that $x \in U(x) \subseteq {\rm pcl} (U(x)) \subseteq
V_{\alpha(x)}$ \cite[Theorem 3.2]{ElD1}. Then, $\{ U(x) \colon
x \in X \}$ is a preopen cover of the p-closed space $X$ and
hence there exists a finite amount of points, say, $x_1, x_2,
\ldots, x_n$ such that $X = \cup_{i=1}^{n} {\rm pcl} (U(x_{i}))
= \cup_{i=1}^{n} V_{\alpha(x_{i})}$. This shows that $X$ is
compact. $\Box$

\section{p-closed subspaces}\label{s3}

Recall that a topological space $(X,\tau)$ is called {\em
hyperconnected} if every open subset of $X$ is dense. In the
opposite case, $X$ is called {\em hyperdisconnected}. A set $A$
is called {\em semi-regular} \cite{DMN1} if it is both semi-open
and semi-closed. Di Maio and Noiri \cite{DMN1} have shown that
a set $A$ is semi-regular if and only if there exists a regular
open set $U$ with $U \subseteq A \subseteq {\rm cl} (U)$. Cameron
\cite{C1} used the term regular semi-open for a semi-regular set.

\begin{lemma}\label{l2a}
{\rm \cite[Mashhour et al.]{MM2}} Let $A$ and $B$ be subsets of
a topological space $(X,\tau)$. 

(1) If $A \in PO(X)$ and $B \in SO(X)$, then $A \cap B \in
PO(B)$.

(2) If $A \in PO(B)$ and $B \in PO(X)$, then $A \in PO(X)$.
\end{lemma}

\begin{lemma}\label{l2}
Let $B \subseteq A \subseteq X$ and $A \in SO(X)$. Then, ${\rm
pcl}_{A} (B) \subseteq {\rm pcl}_{X} (B)$.
\end{lemma}

\begin{theorem}\label{t5}
If every proper semi-regular subspace of a hyperdisconnected
topological space $(X,\tau)$ is p-closed, then $X$ is also
p-closed.
\end{theorem}

{\em Proof.} Since $(X,\tau)$ is not hyperconnected, then there
exists a proper semi-regular set $A$. Let $\{ A_i \}_{i \in I}$
be any preopen cover of $X$. Since $A$ is semi-open, then by
Lemma~\ref{l2a} $A_i \cap A = B_i \in PO(A,\tau|A)$. Then $\{ B_i
\}_{i \in I}$ is a preopen cover of the p-closed space
$(A,\tau|A)$. Then, there exists a finite subset $F$ of $I$ such
that $A = \cup_{i \in F} {\rm pcl}_{A} (B_{i}) \subseteq \cup_{i
\in F} {\rm pcl} (B_{i})$ (by Lemma~\ref{l2}). Therefore, we have
$A \subseteq \cup_{i \in F} {\rm pcl} (A_{i})$. Since $A$ is
semi-regular, $X \setminus A$ is also semi-regular and by a
similar argument we can find a finite subset $G$ of $I$ such that
$X \setminus A \subseteq \cup_{i \in G} {\rm pcl} (A_{i})$.
Hence, $X = \cup_{i \in F \cup G} {\rm pcl} (A_{i})$. This shows
that $X$ is p-closed. $\Box$

\begin{theorem}\label{t6}
If there exists a proper semi-regular subset $A$ of a topological
space $(X,\tau)$ such that $A$ and $X \setminus A$ are
p-closed subspaces, then $X$ is also p-closed.
\end{theorem}

{\em Proof.} Similar to the one of Theorem~\ref{t5} and hence
omitted. $\Box$

\begin{lemma}\label{l3}
Let $A \subseteq B \subseteq X$ and $B \in PO(X)$. If $A \in
PO(B)$, then ${\rm pcl} (A) \subseteq {\rm pcl}_{B} (A)$.
\end{lemma}

\begin{theorem}\label{t7}
If $(X,\tau)$ is a p-closed spaces and $A$ is preregular (i.e.\
both preopen and preclosed), then $(A,\tau|A)$ is also p-closed
(as a subspace).
\end{theorem}

{\em Proof.} Let $\{ A_i \}_{i \in I}$ be any preopen cover of
$(A,\tau|A)$. By Lemma~\ref{l3}, $A_i in PO(X)$ for each $i \in
I$ and $\{ A_i \colon i \in I \} \cup (X \setminus A) = X$. Since
$X$ is p-closed, there exists a finite subset $F$ of $I$ such
that $X = (\cup_{i \in F} {\rm pcl}_{X} (A_{i})) \cup (X
\setminus A$; hence $A = \cup_{i \in F} {\rm pcl}_{X} (A_{i})$.
For each $i \in F$, we have by Lemma~\ref{l3}, ${\rm pcl}_{X}
(A_{i}) \subseteq {\rm pcl}_{A} (A_{i})$ and $A = \cup_{i \in F}
{\rm pcl}_{A} (A_{i})$. Therefore, $A$ is a p-closed subspace.
$\Box$

\begin{example}
{\em An open, even a $\delta$-open subset of a p-closed space
need not be p-closed (as a subspace). Consider any infinite set
$X$ with the point excluded topology. Since the only preopen set
containing the excluded point is the whole space $X$, then the
space in question is p-closed. However, the (infinite) set of
isolated points of $X$ is not p-closed.}
\end{example}

\section{Sets which are p-closed relative to a space}\label{s4}

A subset $S$ of a topological space $(X,\tau)$ is said to be {\em
p-closed relative to $X$} if for every open cover $\{ V_{\alpha}
\colon \alpha \in A \}$ of $S$ by preopen subsets of $(X,\tau)$,
there exists a finite subset $A_0$ of $A$ such that $X = \cup \{
{\rm pcl} (V_{\alpha}) \colon \alpha \in A_0 \}$.

\begin{theorem}\label{t43}
For a topological space $(X,\tau)$ the following conditions are
equivalent:

{\rm (a)} $S$ is p-closed relative to $X$,

{\rm (b)} every maximal filter base on $X$ which meets $S$ \ptc
s to some point of $S$,

{\rm (c)} every filter base on $X$ which meets $S$ \pta s at some
point of $S$,

{\rm (d)} for every family $\{ V_{\alpha} \colon \alpha \in A \}$
of preclosed subsets of $(X,\tau)$ such that $[\cap \{ V_{\alpha}
\colon \alpha \in A \}] \cap S = \emptyset$, there exists a
finite subset $A_0$ of $A$ such that $[\cap \{ {\rm pint}
(V_{\alpha}) \colon \alpha \in A_0 \}] \cap S = \emptyset$.
\end{theorem}

A point $x \in X$ is said to be a {\em pre-$\theta$-accumulation
point} of a subset $A$ of a topological space $(X,\tau)$ if ${\rm
pcl} (U) \cap A \not= \emptyset$ for every $U \in PO(X,x)$. The
set of all pre-$\theta$-accumulation points of $A$ is called the
{\em pre-$\theta$-closure} of $A$ and is denoted by ${\rm
pcl}_{\theta} (A)$. A subset $A$ of a topological space
$(X,\tau)$ is said to be {\em pre-$\theta$-closed} if ${\rm
pcl}_{\theta} (A) = A$. The complement of a pre-$\theta$-closed
set is called {\em pre-$\theta$-open}.

\begin{proposition}\label{p41}
Let $A$ be a subset $A$ of a topological space $(X,\tau)$. Then:

{\rm (i)} If $A \in PO(X)$, then ${\rm pcl} (A) = {\rm
pcl}_{\theta} (A)$.

{\rm (ii)} If $A$ is preregular, then $A$ is pre-$\theta$-closed.

{\rm (iii)} If $A \in SO(X)$, then ${\rm pcl} (A) = {\rm cl}
(A)$.
\end{proposition}

\begin{theorem}\label{tn1}
If $X$ is a p-closed space, then every pre-$\theta$-open cover
of $X$ has a finite subcover.
\end{theorem}

{\em Proof.} Let $\{ V_{\alpha} \colon \alpha \in A \}$ be any
cover of $X$ by pre-$\theta$-open subsets of $X$. For each $x \in
X$, there exists $\alpha{(x)} \in A$ such that $x \in
V_{\alpha{(x)}}$. Since $V_{\alpha(x)}$ is pre-$\theta$-open,
there exists $V_x \in PO(X)$ such that $x \in V_x \subseteq {\rm
pcl} (V_{x}) \subseteq V_{\alpha(x)}$. The family $\{ V_x \colon
x \in X \}$ is a preopen cover of $X$. Since $X$ is p-closed,
there exists a finite number of points, say, $x_1,x_2, \ldots,
x_n$ in $X$ such that $X = \cup \{ {\rm pcl} (V_{x_i}) \colon i =
1,2,\ldots, n \}$. Therefore, we obtain that $X = \cup \{
V_{\alpha(x_i)} \colon i = 1,2,\ldots,n \}$. $\Box$

{\bf Question.} Is the converse in Theorem~\ref{tn1} true?

\begin{theorem}\label{tn2}
Let $A,B$ be subsets of a space $X$. If $A$ is
pre-$\theta$-closed and $B$ is p-closed relative to $X$, then $A
\cap B$ is p-closed relative to $X$.
\end{theorem}

{\em Proof.} Let $\{ V_{\alpha} \colon \alpha \in A \}$ be any
cover of $A \cap B$ by preopen subsets of $X$. Since $X \setminus
A$ is pre-$\theta$-open, for each $x \in B \setminus A$ there
exists $W_x \in PO(X,x)$ such that ${\rm pcl} (W_x) \subseteq X
\setminus A$. The family $\{ W_x \colon x \in B \setminus A \}
\cup \{ V_{\alpha} \colon \alpha \in A \}$ is a cover of $B$ by
preopen sets of $X$. Since $B$ is p-closed relative to $X$, there
exist a finite number of points, say, $x_1,x_2,\ldots,x_n$ in $B
\setminus A$ and a finite subset $A_0$ of $A$ such that $$B
\subseteq [ \cup_{i=1}^{n} {\rm pcl} (W_{x_{i}}) ] \cup [
\cup_{\alpha \in A_{0}} {\rm pcl} (V_{\alpha}) ].$$ Since ${\rm
pcl} (W_{x_{i}}) \cap A = \emptyset$ for each $i$, we obtain that
$A \cap B \subseteq \cup \{ {\rm pcl} (V_{\alpha}) \colon \alpha
\in A_0 \}$. This shows that $A \cap B$ is p-closed relative to
$X$. $\Box$

\begin{corollary}\label{c45}
If $K$ is pre-$\theta$-closed set of a p-closed space $(X,\tau)$,
then $K$ is p-closed relative to $X$.
\end{corollary}

{\bf Question.} If in a topological space $(X,\tau)$ every proper
pre-$\theta$-closed set is p-closed relative to $X$, is $X$
necessarily p-closed?

A topological space $(X,\tau)$ is called {\em preconnected}
\cite{P1} if $X$ can not be expressed as the union of two
disjoint preopen sets. In the opposite case, $X$ is called {\em
predisconnected}. Note that every preconnected space is
irresolvable but not vice versa.

\begin{theorem}\label{tn3}
Let $X$ be a predisconnected space. Then $X$ is p-closed if and
only if every preregular subset of $X$ is p-closed relative to
$X$.
\end{theorem}

{\em Proof.} {\em Necessity.} Every preregular set is
pre-$\theta$-closed by Proposition~\ref{p41}. Since $X$ is
p-closed, the proof is completed by Corollary~\ref{c45}.

{\em Sufficiency.} Let $\{ V_{\alpha} \colon \alpha \in A \}$ be
a preopen cover of $X$. Since $X$ is predisconnected, there
exists a proper preregular subset $A$ of $X$. By our hypothesis,
$A$ and $X \setminus A$ are p-closed relative to $X$. There exist
finite subsets $A_1$ and $A_2$ of $A$ such that $$A \subseteq
\cup_{\alpha \in A_1} {\rm pcl} (V_{\alpha}) \ {\rm and} \ X
\setminus A \subseteq \cup_{\alpha \in A_2} {\rm pcl}
(V_{\alpha}).$$ Therefore, we obtain that $X = \cup \{ {\rm pcl}
(V_{\alpha} \colon \alpha \in A_1 \cup A_2 \}$. $\Box$

\begin{theorem}\label{tn4}
If there exists a proper preregular subset $A$ of a topological
space $(X,\tau)$ such that $A$ and $X \setminus A$ are p-closed
relative to $X$, then $X$ is p-closed.
\end{theorem}

{\em Proof.} This proof is similar to the one of
Theorem~\ref{tn3} and hence omitted. $\Box$

\begin{theorem}\label{tn5}
Let $X_0$ be a semi-open subset of a topological space
$(X,\tau)$. If $X_0$ is a p-closed space, then it is p-closed
relative to $X$.
\end{theorem}

{\em Proof.} Let $\{ V_{\alpha} \colon \alpha \in A \}$ be any
cover of $X_0$ by preopen subsets of $X$. Since $X_0 \in SO(X)$,
by Lemma~\ref{l2a} we have that $X_0 \cap V_{\alpha} = W_{\alpha}
\in PO(X_0)$ for each $\alpha \in A$. Therefore, $\{ W_{\alpha}
\colon \alpha \in A \}$ is a preopen cover of $X_0$. since $X_0$
is p-closed, there exists a finite subset $A_0$ of $A$ such that
$X_0 = \cup \{ {\rm pcl}_{X_0} (W_{\alpha}) \colon \alpha \in A_0
\}$. By Lemma~\ref{l2}, we obtain that $X_0 \subseteq \cup \{
{\rm pcl} (W_{\alpha}) \colon \alpha \in A_0 \} \subseteq \cup
\{ {\rm pcl} (V_{\alpha}) \colon \alpha \in A_0 \}$. This shows
that $X_0$ is p-closed relative to $X$. $\Box$

\begin{theorem}\label{tn6}
Let $X_0$ be a preopen subset of a topological space $(X,\tau)$.
If $X_0$ is a p-closed relative to $X$, then it is a p-closed
subspace of $X$.
\end{theorem}

{\em Proof.} Let $\{ V_{\alpha} \colon \alpha \in A \}$ be any
cover of $X_0$ by preopen subsets of $X$. Since $X_0 \in PO(X)$,
by Lemma~\ref{l2a}, $V_{\alpha} \in PO(X)$ for each $\alpha \in
A$. Since $X_0$ is p-closed relative to $X$, there exists a
finite subset $A_0$ of $A$ such that $X_0 \subseteq \cup \{ {\rm
pcl} (V_{\alpha}) \colon \alpha \in A_0$. This shows that $X_0$
is p-closed subspace of $X$. $\Box$

\begin{corollary}
Let $X_0$ be an ($\alpha$-)open subset of a topological space
$(X,\tau)$. Then $X_0$ is a p-closed subspace of $X$ if and only
if it is p-closed relative to $X$.

\end{corollary}

{\em Proof.} This is an immediate consequence of
Theorem~\ref{tn5} and Theorem~\ref{tn6}. $\Box$

Recall that a function \fxy\ is called {\em preirresolute}
\cite{RV1} (resp.\ {\em precontinuous} \cite{MAH1}) if $f^{-1}
(V)$ is preopen in $X$ for every preopen (resp.\ open) subsets
$V$ of $Y$.

\begin{lemma}\label{lp1}
{\em \cite{P2}} A function \fxy\ is preirresolute (resp.\
precontinuous) if and only if for each subset $A$ of $X$, $f
({\rm pcl} (A)) \subseteq {\rm pcl} (f(A))$ (resp.\ $f ({\rm pcl}
(A)) \subseteq {\rm cl} (f(A))$).
\end{lemma}

\begin{theorem}
If a function \fxy\ is a preirresolute (resp.\ precontinuous)
surjection and $K$ is p-closed relative to $X$, then $Y$ is
p-closed (resp.\ QHC) relative to $Y$.
\end{theorem}

{\em Proof.} Let $\{ V_{\alpha} \colon \alpha \in A \}$ be any
cover of $f(K)$ by preopen (resp.\ open) subsets of $Y$. Since
$f$ is preirresolute (resp.\ precontinuous), $\{ f^{-1}
(V_{\alpha}) \colon \alpha \in A \}$ is a cover of $K$ by preopen
subsets of $X$, where $K$ is p-closed relative to $X$. Therefore,
there exists a finite subset $A_0$ of $A$ such that $K \subseteq
\cup_{\alpha \in A_0} {\rm pcl} (f^{-1}(V_{\alpha}))$. Since
$f$ is preirresolute (resp.\ precontinuous) and surjective, by
Lemma~\ref{lp1}, we have $$f(K) \subseteq \cup_{\alpha \in A_0}
f({\rm pcl} (f^{-1}(V_{\alpha}))) \subseteq \cup_{\alpha \in A_0}
{\rm pcl} (V_{\alpha})$$ $$({\rm resp.} \ f(K) \subseteq
\cup_{\alpha \in A_0} f({\rm pcl} (f^{-1}(V_{\alpha}))) \subseteq
\cup_{\alpha \in A_0} {\rm cl} (V_{\alpha})).$$

\begin{corollary}
If a function \fxy\ is a preirresolute (resp.\ continuous)
surjection and $X$ is p-closed, then $Y$ is p-closed (resp.\
QHC).
\end{corollary}

\begin{corollary}
{\rm (i)} The property ``p-closed" is topological.

{\rm (ii)} If the product space $\prod_{\alpha \in A} X_{\alpha}$
is p-closed, then $X_{\alpha}$ is p-closed for each $\alpha \in
A$.
\end{corollary}

\begin{remark}
{\em Even finite product of p-closed spaces need not be p-closed;
for consider the product of the space from Example~\ref{e1} (i)
with any two point indiscrete space. This product space shows
that Theorem 3.4.3 from \cite{A1} is wrong, i.e., every proper
preregular subset might be p-closed relative to the space and
still the space might fail to be p-closed. Additionally, Example
3.4.1 from \cite{A1} is also false.}
\end{remark}

\
\begin{center}
Department of Mathematics\\University of Helsinki\\PL 4,
Yliopistonkatu 15\\00014 Helsinki\\Finland\\e-mail: {\tt
dontchev@cc.helsinki.fi}, {\tt dontchev@e-math.ams.org}
\end{center}
\
\begin{center}
Department of Mathematics\\Graz University of
Technology\\Steyrergasse 30\\A-8010 Graz Austria\\e-mail: {\tt
ganster@weyl.math.tu-graz.ac.at}
\end{center}
\
\begin{center}
Department of Mathematics\\Yatsushiro College of Technology\\2627
Hirayama shinmachi\\Yatsushiro-shi\\Kumamoto-ken\\866
Japan\\e-mail: {\tt noiri@as.yatsushiro-nct.ac.jp}
\end{center}
\
\end{document}